\documentclass{amsart}

\theoremstyle{definition}

\theoremstyle{remark}

\numberwithin{equation}{section}

\begin{document}

\title{von Neumann entropy and relative position between subalgebras }

\author{Marie Choda}

\footnote{{\it  2010 Mathematics Subject Classification} :  46L55;  46L37, 46L40}
\keywords{entropy, unitary, density matrix }
\thanks{The author was supported in part by JSPS Grant No.20540209.}

\maketitle

\centerline{\small Osaka Kyoiku University, 
Asahigaoka, Kashiwara 582-8582, Japan} 
\centerline{marie@cc.osaka-kyoiku.ac.jp}



\date{}

\dedicatory{}

\begin{abstract} 
We give a numerical characterization of mutual orthogonality  (that is,  complementarity) 
for subalgebras. 
In order to give such a characterization for mutually orthogonal  subalgebras $A$ and $B$ of 
the $k \times k$ matrix algebra $M_k(\mathbb{C})$, where $A$ and $B$  are  isomorphic to some 
$M_n(\mathbb{C})$ $(n \leq k)$,   
we  consider  a density matrix which is induced from the pair $\{A, B\}$. 
We show that $A$ and $B$ are  mutually orthogonal if and only if the von Neumann entropy of  
the density matrix is  the maximum value $2\log n$,  
which is the logarithm of the dimension of the subfactors. 
\end{abstract}
\smallskip

\section{Introduction} 
There are several notions which describe some relative position between two subalgebras of operator algebras. 
As one of such  notions  between two subalgebras of finite von Neumann 
algebras, Popa introduced the notion of {\it mutually orthogonal subalgebras } (definition below) in \cite{Po}. 
By the terminology {\it complementarity}, 
the same notion is investigated in the theory of quantum systems (see  \cite{Pe2} for example).  

The most primary interest would be the case where two subalgebras of some full matrix algebra, both of  
which  are either maximal abelian or isomorphic to also some full matrix algebra. 
In such the  cases, two subalgebras are connected by a unitary. 

Our motivation for this work arises from the following fact:   

In the previous paper \cite{Ch1}, we defined a constant $h(A|B)$ 
for two subalgebras $A$ and $B$ of  a finite von Neumann algebra, 
and showed  the relative position between maximal abelian subalgebras $A$ and $B$ of $M_n(\mathbb{C})$ 
by using the values of  $h(A|B)$. 
This  $h(A|B)$ is  a slight modification of  Connes-St\o rmer  relative entropy $H(A|B)$ in \cite{CS} (cf. \cite{NS}). 
If $A_1$ and $A_2$ are maximal abelian subalgebras of $M_n(\mathbb{C}),$ 
then there exists a unitary $u$ such that $A_2 = uA_1u^*,$ and then 
$h(A_1| A_2)$ coincides with the entropy $H(b(u))$  defined in \cite{ZSKS} of the unistochastic matrix $b(u)$ 
induced by the unitary $u$. 
As a consequence, we showed that $A_1$ and $A_2$ are  mutually orthogonal if and if $h(A_1| A_2) = H(b(u)) = \log n.$ 
This means that $A_1$ and $A_2$ are  mutually orthogonal if and if $h(A_1| A_2)$ is maximal and  
equals to the  logarithm of the dimension of the subalgebras. 
Related results in the case of subfactors of the type II$_1$ factors are obtained in \cite {Ch2}.  
Here, it does not hold in general that  $H(A_1| A_2) = H(b(u))$ 
(see, for example \cite [Appendix] {PSW} ). 

On the other hand, Petz showed in \cite{Pe2} for subalgebras $A$ and $B$ of $M_n(\mathbb{C})$ 
that if $A$ is homogeneous and abelian, then  $H(A|B)$  is 
maximal if and only if $A$ and $B$ are complementary. 
Here  homogeneous means that all minimal projections of $A$ have the same trace. 
Also he remarked  that  Connes-St\o rmer relative entropy cannot characterize the complementarity of 
subalgebras in the general case. 
\smallskip 

In this paper, we study the case when the subalgebras  $A$ and $B$ in question are  isomorphic to 
some $M_n(\mathbb{C})$. 
We introduce  some density matrix arising from  the pair $\{A, B\},$ and 
we show that the von Neumann entropy of the density matrix gives  a  characterization of 
the mutual orthogonality (that is the complementarity). 
\smallskip 

In order  to define the  entropy for automorphisms of operator algebras, 
two kind of notion of {\it a finite partition of unity} played an important role. 
One was used by Connes-St\o rmer, and it corresponds to a  finite measurable partition of 
a given space in the ergodic theory (see \cite {NS} \cite {OP} for example). 
The other was used by Alicci and Fannes in \cite{AF} and it is called  a {\it finite operational 
partition of unity}. 
Here, we apply the latter, that is operational partition of unity, 
and we give a numerical characterization for pairs of mutually orthogonal 
 subalgebras which are both isomorphic to some full matrix algebra of the same size. 
\smallskip 

The paper is organized as follows. After preliminaries on 
basic notions in Section 2, in Section 3  we define  some density matrix 
which is closely related to subfactors $A$ and $B$  which are both isomorphic to some 
$M_n(\mathbb{C})$, and  we   show that $A$ and $B$ are  mutually orthogonal 
if and only if the von Neumann entropy of  
the density matrix is the maximum value $2\log n$,  
which is the logarithm of the dimension of the subfactors. 
\vskip 0.3cm

\section{Preliminaries}
Here we summarize notations, terminologies and 
basic facts. 

Let $M$ be a finite von  Neumann algebra acting on a separable Hilbert space,  
and $\tau$ be a fixed normal faithful tracial state. 
In the case where $M$ is the algebra $M_n(\mathbb{C})$ of $n \times n$ matrices, 
$\tau (x) = {\rm Tr}(x)/ n, $ where  {\rm Tr} the usual standard trace on $M_n(\mathbb{C})$. 
The norm $\Vert x \Vert_\tau$ is given by $\Vert x \Vert_\tau = \tau(x^*x)^{1/2}$ for all $x \in M$.  
By a von Neumann subalgebra $A$ of $M,$  we mean that $A$ is a $*$-subalgebra closed 
in the weak operator topology, the unit of which is the same with the unit of $M$.  
A conditional expectation of $M$  onto 
a von Neumann subalgebra  $A$ of $M$ is  
a completely positive linear map $E_A : M \to A$ with 
$E_A(axb) = aE_A(x)b$ for all $x \in M$ and $a, b \in A$.  
In the case of a finite von Neumann algebra $M$ with 
a  faithful normal tracial state $\tau$, 
there exists always a unique faithful normal conditional expectation $E_A$ 
of $M$ onto a von Neumann subalgebra  $A$ of $M$ 
such that $\tau(xa) = \tau(E_A(x)a)$ for all $x \in M$ and $a \in A$. 
It is called the conditional expectation with respect to $\tau$.  

\subsection{\bf Mutually orthogonal (or complementary) subalgebras.} 
Let $A$ and $B$ be von Neumann subalgebras of $M$. 
In \cite [Lemma 2.1] {Po}, Popa showed that the following conditions are equivalent.  
\begin{enumerate}
\item$\tau(ab) = 0$ for $a \in A, b \in B$ with $\tau(a) = \tau(b) = 0$;
\item $\tau(ab) = \tau(a) \tau(b)$ for all $a \in A, b \in B$;
\item $\Vert ab \Vert_{\tau} = \Vert a \Vert_{\tau} \Vert b \Vert_{\tau} $ for all $a \in A, b \in B$;
\item $E_A E_B(x)  = \tau(x) 1_M,$ for all $x \in M$; 
\item $E_A (B) \subset \mathbb{C} 1_M$.  
\end{enumerate}
Moreover (1) - (5) are equivalent with the analogue conditions obtained by interchanging $A$ with $B$. 

\smallskip
Two  von Neumann subalgebras $A$ and $B$ of $M$ are called {\it mutually orthogonal} 
if one of the above conditions (1) - (5) is satisfied (\cite [Definition 2.2] {Po}). 

Mutually orthogonal subalgebras are also called  {\it complementary  subalgebras}, 
(see, \cite{Pe1} \cite{Pe2} for example). 

\subsection{\bf Density matrix and von Neumann entropy.} 
By a density matrix, we mean a positive semidefinite matrix $\rho$ 
such that ${\rm Tr}(\rho) = 1$. 
To a density matrix  $\rho$, the von Neumann entropy $S(\rho)$ is given by 
$S(\rho) = {\rm Tr}(\eta(\rho))$. Here, $\eta$ is   defined on the interval $[0,1]$ by 
$$\eta(t) = -t \log t \quad  (0 < t \leq 1) \quad \text{and} \quad \eta(0) = 0.$$ 
\smallskip 

\section{Main results} 
Let $ M_n(\mathbb{C}) $  be the algebra of $n\times n$ complex  matrices,  and let 
${\rm Tr}$ be the trace of $ M_n(\mathbb{C}) $ with ${\rm Tr}(p) = 1 $ for every minimal projection $p$. 
Let $L$ be a finite von Neumann algebra, and 
 let $\tau_L$ be a fixed normal  faithful tracial state. 
 
We let $M = M_n(\mathbb{C}) \otimes L,$ and let $\tau_M = {\rm Tr}/ n \otimes \tau_L$. 
\smallskip

\subsection{} We consider the subalgebra $N = M_n(\mathbb{C}) \otimes 1_L$ of $M$. 
In this case, the conditional expectation $E_N$ with respect to $\tau_M$ satisfies that 
$$E_N(x \otimes y) = \tau_L(y)x \otimes 1_L, \quad x \in M_n(\mathbb{C}), \quad y \in L.$$

The following lemma is an easy consequence from the definition, and it is essential to our study. 

\subsubsection{\bf Lemma} 
{\it 
Let $N = M_n(\mathbb{C}) \otimes 1_L$ and let $u \in M$ be a unitary operator. Then 
$N$ and $uNu^*$ are mutually orthogonal if and only if 
$$E_N(u^*(a \otimes 1_L)u) = \tau_M(a \otimes 1_L) 1_M 
 = \frac{{\rm Tr}(a)} n 1_M,\quad \text{for all} \quad a \in  M_n(\mathbb{C}).$$} 

\begin{proof} 
Assume that $N$ and $uNu^*$ are mutually orthogonal, that is, 
$$E_NE_{uNu^*}(x) = E_{uNu^*}E_N(x) = \tau_M(x) 1_M, \quad \text{for all} \quad x \in M.$$ 
Then 
$uE_N(u^*xu)u^* = E_{uNu^*}(E_N(x)) =  \tau_M(x) 1_M$, for all $x \in N$. 
This implies that 
$$E_N(u^*xu) =  \tau_M(x) 1_M, \quad \text{for \ all} \quad x \in N.$$

Conversely, 
assume that $E_N(u^*xu) =  \tau_M(x) 1_M$, for all $x \in N$. Then 
$$E_{uNu^*}(x) = uE_N(u^*xu)u^* = \tau_M(x) 1_M$$ 
for all $x \in N$. 
Hence 
$$E_{uNu^*}E_N(x) = \tau_M(x) 1_M \quad \text{ for \ all} \ x \in M$$ 
so that $N$ and $uNu^*$ are mutually orthogonal.
\end{proof}
\smallskip

\subsection{}
Let $\{e_{ij} ; i,j = 1, \cdots, n \}$ be a system of matrix units of $M_n(\mathbb{C}),$ 
so that 
$$e_{ij}^* = e_{ji}, \quad e_{ij} e_{st} = \delta_{js} e_{it}, \quad \sum_{i = 1}^n e_{ii} = 1_{M_n(\mathbb{C})}.$$ 
Then each  $x$ in $M = M_n(\mathbb{C}) \otimes L$ is written in the unique form: 
$$x = \sum_{i,j = 1}^n e_{ij} \otimes x_{ij}, \quad x_{ij} \in L,$$
and  $u = \sum_{i,j = 1}^n e_{ij} \otimes u_{ij}$  is a unitary in $M$ if and only if 
$$\sum_{j=1}^n u_{ij} u_{kj}^* = \delta_{ik} 1_{L}\quad \text{and} \quad 
\sum_{i=1}^n u_{ij}^* u_{ik} = \delta_{jk} 1_L . $$
\smallskip

We give a characterization for a unitary $u \in M$ to satisfy that $N$ and $uNu^*$ are mutually orthogonal. 
\smallskip

\subsubsection{\bf Theorem.}  
{\it  
Assume that a von Neumann subalgebra $N$ of $M$ is given by  $N =  M_n(\mathbb{C}) \otimes 1_L$ 
and let $ u \in M$ be  unitary.  
Then $N$ and $uNu^*$ are mutually orthogonal if and only if 
$$\tau_L(u_{ij}^* u_{kl}) = \delta_{ik}\delta_{jl} \frac 1n, \quad \text{for  all} \quad i,j,k,l = 1, \cdots, n.$$
}

\begin{proof}  
Assume that  $N$ and $uNu^*$ are mutually orthogonal. 
Then by Lemma 3.1.1 
$$ E_N(u^*(e_{ij} \otimes 1_L) u ) = \delta_{ij}\frac 1n 1_M.$$ 
On the other hand, since
$$u^*(e_{ij} \otimes 1_L) u = \sum_{l, t = 1}^n e_{lt} \otimes u_{il}^* u_{jt}, \ \text{for \ all} \ i, j  = 1, \cdots, n,$$ 
by applying  that $E_N(x \otimes y) = \tau_L(y)x \otimes 1_L,$ 
we have that 
$$ E_N(u^*(e_{ij} \otimes 1_L) u )
 = \sum_{l, t = 1}^n \tau_L( u_{il}^* u_{jt}) e_{lt} \otimes 1_L.$$
Hence  
$\sum_{l, t = 1}^n \tau_L( u_{il}^* u_{jt}) e_{lt} = \delta_{ij}\frac 1n 1_{M_n(\mathbb{C})}$. 
This means that 
$$\tau_L(u_{ij}^* u_{kl}) = \delta_{ik}\delta_{jl} \frac 1n \quad \text{for  all} \quad i,j,k,l = 1, \cdots, n.$$

Conversely, assume that $\tau_L(u_{ij}^* u_{kl}) = \delta_{ik}\delta_{jl} \frac 1n$ for all $i,j,k,l = 1, \cdots, n$. 
Then we have that 
$${E_N (u^*(e_{ij} \otimes 1_L)u) } =  \sum_{l, t = 1}^n e_{lt} \otimes \tau_L(u_{il}^*u_{jt}) 1_L 
 =  \sum_{l = 1}^n e_{ll} \otimes \delta_{ij} \frac 1n 1_L \nonumber = \delta_{ij} \frac 1n 1_M
$$
for all $i,j,k,l = 1, \cdots, n.$ 
Hence   $N$ and $uNu^*$ are mutually orthogonal, by Lemma 3.1.1. 
\end{proof}

\subsubsection{\bf Note}  
Theorem 3.2.1 implies that if $N =  M_n(\mathbb{C}) \otimes 1_L$  and if $N$ and $uNu^*$ are mutually orthogonal 
for some unitary $u \in M = M_n(\mathbb{C}) \otimes L,$ then 
the set $\{  u_{ij}  / {\sqrt n} \ ; i,j = 1, \cdots, n\} \subset L$ has to be an orthonormal system with respect to 
the inner product induced by $\tau_L$ so that $\dim(L) \geq n^2$.

\subsection{ Entropy associated to an inner conjugate pair of subfactors} 
In order to give a numerical characterization for mutually orthogonal subalgebras which are all isomorphic to 
$M_n(\mathbb{C})$, 
we  apply the notion of a  finite operational partition $X$ of unity of size $k$  and the density matrix $\rho_\phi[X]$ 
which were introduced by Alicki and Fannes in \cite {AF}. 

\subsubsection{\bf Finite  operational partition} 
Let $A$ be a unital $C^*$-algebra. A {\it finite  operational partition of unity of size $k$ } is a set 
$X = \{x_1, ..., x_k \}$ of elements of $A$ satisfying 
$$\sum_i^k x_i^* x_i = 1_A.$$

We remark that a similar terminology  "finite partition" is usually used  in the different following form: 
A finite subset $ \{x_1, ..., x_k \}$ in  $A$ is called  a finite partition of unity if 
they are nonnegative operators in $A$ such that $1_A = \sum_{i=1}^n  x_i$. See \cite{NS} or \cite{OP}. 

\subsubsection{{\bf Density matrix} $\rho[X]$} 

Let $\phi$ be a state of $A$. 
To a finite operational partition $X$  of unity of size $k$,   
we associate a $k \times k$ density matrix  $\rho_\phi[X]$ such that 
the $(i,j)$-coefficient $\rho_\phi[X] (i,j)$ of $\rho_\phi[X]$ is given by 
$$ \rho_\phi[X] (i,j) = \phi(x_j^*x_i), \quad i,j = 1, \cdots, k.$$ 
In the case that $A$ is a finite von Neumann algebra and that $\phi$ is a given tracial state $\tau$ of $A,$  
then we denote $\rho_\tau[X]$ simply by $\rho[X].$
\smallskip

\subsubsection{{\bf Finite operational partition  induced by a unitary} $u$} 
Now let $ M_n(\mathbb{C}) $  be the algebra of $n\times n$ complex  matrices  and let 
${\rm Tr}$ be the trace with ${\rm Tr}(p) = 1 $ for every minimal projection $p$. 
Let $L$ be a finite von Neumann algebra, and
 let $\tau_L$ be a fixed normal  faithful tracial state. 
Let $M = M_n(\mathbb{C}) \otimes L,$ and let $\tau_M = {\rm Tr}/ n \otimes \tau_L$. 
Let $u$ be a unitary in $M_n(\mathbb{C}) \otimes L,$ 
and let $u = \sum_{i,j} e_{ij} \otimes u_{ij}, \ (u_{ij} \in L),$ 
where $\{e_{ij}\}_{i,j = 1, \cdots, n}$ is a set of matrix units of $M_n(\mathbb{C})$. 
We consider the set 
$$U = \{\frac 1{\sqrt n} {u_{ij}} \ ; \ i, j = 1, \cdots, n\}.$$ 
It is not so essential, but 
we renumber the elements of $U$ for the sake of convenience. 
For example, if $kn+1 \leq i \leq (k+1)n, $ for some $k = 0, 1, \cdots, n-1,$ then we put 
$$u_i = \frac 1{\sqrt n} {u_{i-kn \ k+1}}.$$
It is clear the correspondence $i \longleftrightarrow (i-kn,  k+1)$ for some $k = 0, 1, \cdots, n-1$ 
is one to one. 
Since $u$ is a unitary, clearly  the set $U$ is a finite  operational partition  of unity of size $n^2$. 
We call this set $U$ the  {\it finite operational partition  of unity induced by} $u$. 
\smallskip

\subsubsection{{\bf von Neumann entropy} $S(\rho[U] )$} 
We consider 
the von Neumann entropy $S(\rho_\phi[U] )$ of the density operator $\rho_\phi[U] $ 
in order to characterize the mutual orthogonality  for subfactors. 
So,  we assume that our state $\phi$ is the  given normalized trace and 
$$ S(\rho[U] ) = {\rm Tr}(\eta(\rho[U] )).$$  
\smallskip

\subsubsection{\bf Theorem.}  
{\it 
Let $L$ be a finite von Neumann algebra and let $\tau_L$ be a normalized trace of  $L$. 
We let  $ M =  M_n(\mathbb{C}) \otimes L$ and  $\tau = {\rm Tr} / n  \otimes \tau_L$.
Assume that $N =  M_n(\mathbb{C}) \otimes 1_L$   
and that $u$ is a unitary operator in $M$.  
Then  the following conditions are equivalent: 
\begin{enumerate}
 \item $N$ and $uNu^*$ are mutually orthogonal; 
 \item $n^2 \rho[U] $ is the $n^2 \times n^2$ identity matrix; 
 \item $ S(\rho[U] ) = 2 \log n = \log\dim  N. $  
\end{enumerate}
Here $U$ is the  finite operational partition  of unity induced by $u$. 
}
\smallskip

\begin{proof}  
First we remark that 
$$\rho[U](i,j) = \frac 1n \tau(u_{j-ln, \ l+1}^* u_{i-kn, \  k+1})$$ 
where  $u_i = (1 / {\sqrt n}) u_{i-kn, \ k+1}  $,  
for some $k = 0, 1, \cdots, n-1$ with $kn+1 \leq i \leq (k+1)n, $ 
and 
$u_j = ( 1 / {\sqrt n}) {u_{j-ln \ l+1}}$,  
for some $l = 0, 1, \cdots, n-1$ with 
$ln+1 \leq j \leq (lk+1)n $. 
\smallskip

(1) $\Rightarrow$ (2):   
Assume that $N$ and $uNu^*$ are mutually orthogonal. 
Then by Theorem 3.2.1 and by the definition of $\rho[U],$ 
the $n^2 \times n^2$ density matrix $\rho[U]$ is the diagonal 
matrix such that 
$$\rho[U] (i,i) = \frac 1{n^2} \quad \text{for} \ i = 1, 2, \cdots, n^2.$$
\smallskip

(2) $\Rightarrow$ (3):    
Clearly, the  the von Neumann entropy $S(\rho [U] ) = 2 \log n$ and it is the dimension of $N$. 
\smallskip

(3) $\Rightarrow$ (2):  
Assume that $S(\rho[U] ) =  \log n^2$. 
Let $(\lambda_1, \cdots, \lambda_{n^2})$ be  an eigenvalue sequence of $\rho[U] $ and let 
$(p_1,  \cdots, p_{n^2})$ be  the corresponding sequence of the minimal projections. 
Then there exists a  $n^2 \times n^2$  unitary  matrix $w$ so that 
$$w\rho[U] w^* = \sum_{i = 1} ^{n^2} \lambda_i p_i.$$ 
Since 
$$ \log n^2 = S(\rho[U] )= \sum_{i = 1} ^{n^2} \eta(\lambda_i ), $$
it implies that, by the concavity of the function $\eta$, 
$$\lambda_i = \frac 1{n^2} \quad \text{for all} \quad i = 1, 2, \cdots, n^2$$
so that 
$$w\rho[U] w^* = \frac 1{n^2} 1_{M_{n^2}(\mathbb{C})}.$$
Hence (2) holds. 
\smallskip

(2) $\Rightarrow$ (1):   
By the definition of $\rho[U]$ and the condition (2), we have that 
$$\delta_{ij} \frac 1{n^2} = \rho[U](i,j) = \frac 1n \tau(u_{j-ln \ l+1}^* u_{i-kn \  k+1}).$$ 

This relation corresponds that   
$\tau_L(u_{ij}^* u_{kl}) = \delta_{ik}\delta_{jl} \frac 1n$. 
Hence by Theorem 3.2.1,   $N$ and $uNu^*$ are mutually orthogonal. 
\end{proof}
\vskip 0.3cm

\subsubsection{\bf Note.}  
Theorem 3.3.5 means that the mutually orthogonality for inner conjugate 
subfactors are characterized by  the maximum value 
$\log \dim$ of the subfactors.  

In fact, since the density matrix $\rho[U]$ is a $n^2 \times n^2$ matrix and  
the function $\eta$ is operator concave, the value $2 \log n$ is the maximum. 

\subsubsection{\bf Note.}  
The proof  shows that the statement of Theorem 3.3.5 does not depend on any choice of a matrix units.

\subsection{\bf Subfactors of matrix algebras}
Let $A$ and $B$ be subalgebras of  $M_k(\mathbb{C})$ and assume that both subalgebras are isomorphic to 
$M_n(\mathbb{C})$. Then $k = mn$. We can assume that 
$M_k(\mathbb{C}) = M_n(\mathbb{C}) \otimes M_m(\mathbb{C})$ and 
$A = M_n(\mathbb{C}) \otimes \mathbb{C}1$.   
There exists a unitary matrix $u \in M_k(\mathbb{C})$ such that  $B = uAu^*$.    
We denote by $u(A,B)$ this unitary  and also by $U(A,B)$  the  finite operational 
partition  of unity induced by $u(A,B)$.  
Then we have the followings:  

\subsubsection{} 
Petz's  characterization of complementarity was given in (\cite [Theorem 4] {Pe1}): 
The subalgebra $u(1 \otimes M_m(\mathbb{C}) ) u^*$ is complementary to $1 \otimes M_m(\mathbb{C})$ 
if and only if 
$$\frac mn \sum_{i,j = 1}^n |u_{ij} >< u_{ij}| = 1.$$ 
When n = m this condition means that $\{u_{ij}\}_{ij}$  is an orthonormal basis in $M_n(\mathbb{C})$ 
with respect to the inner product by ${\rm Tr}$. 
\vskip 0.3cm

\smallskip
Our characterization is the following Corollary of Theorem 3.3.6 by letting $L =  M_m(\mathbb{C})$. 

\subsubsection{\bf Corollary.}
{\it
Let $A$ and $B$ be subalgebras of  $M_{k}(\mathbb{C})$ and assume that both subalgebras are isomorphic to 
$M_n(\mathbb{C})$. Then $A$ and $B$ are mutually orthogonal  if and only if 
$$ S(\rho[U(A,B)] ) = 2 \log n = \log(\dim A). $$ 
}
\smallskip
\subsubsection{\bf Note} 
In the above 3.4.1 and 3.4.2, the numbers $m$ and $n$ should be  $m \geq n.$ 
\smallskip

\subsubsection{\bf Comparison with the case of maximal abelian subalgebras.} 
We remark that  Corollary 3.4.2  corresponds to \cite [Corollary 3.2, Corollary 3.3] {Ch1}: 
\smallskip

Assume that $A$ and $B$ are maximal abelian subalgebras of  $M_{n}(\mathbb{C})$. 
Then there exists  a unitary $u$ in $M_{n}(\mathbb{C})$ with $uAu^* = B$, 
and we have that 
\begin{enumerate}
 \item $h(A \mid B) = H(b(u)).$ 
  \item 
  $A$ and $B$ are mutually orthogonal  if and only if 
    $$h(A \mid B) = \log n = \log(\dim A).$$ 
\end{enumerate}

Here, $h(A \mid B) $ is the conditional relative entropy for $A$ and $B$ in \cite{Ch1}  and 
$ H(b(u)) $ is the entropy  for the unistochastic operator $b(u)$ induced by the unitary $u$  in \cite{ZSKS}. 
\smallskip

This means that $A$ and $B$ are mutually orthogonal  if and only if 
$h(A \mid B)$ takes the maximum value $\log(\dim A)$, 
because $\log n$ is the maximum value by the definition of $H(b(u))$ and by 
the property of the function $\eta$.  
\smallskip

\bibliographystyle{amsplain}

\end{document}